\documentclass{article}[12pt]
\usepackage{natbib}
\usepackage{authblk}
\usepackage{graphicx}
\usepackage{amsmath}
\usepackage[a4paper,top=2cm,bottom=2cm,left=2cm,right=2cm]{geometry}
\usepackage{mathrsfs}
\usepackage{amsthm}
\usepackage{amssymb}
\usepackage{geometry}
\usepackage{setspace}
\usepackage{graphicx}
\usepackage{extarrows}
\usepackage{float}
\usepackage{color}
\usepackage{enumitem}
\usepackage[title]{appendix}
\newcommand{\keywords}[1]{\textbf{\textit{Keywords:}} #1}
\def\du{\ensuremath{\mathrm{d}}}
\numberwithin{equation}{section}

\theoremstyle{plain}
\newtheorem{theorem}{Theorem}
\newtheorem{lemma}{Lemma}

\newtheorem{proposition}{Proposition}
\newtheorem{definition}{Definition}

\usepackage{hyperref}
\usepackage{url}

\begin{document}

\title{Nearly minimax empirical Bayesian prediction of independent Poisson observables}
\author[1,*]{Xiao Li}
\affil[1]{\small Department of Mathematical Informatics, Graduate School of Information Science and Technology,\par The University of Tokyo, 7-3-1 Hongo, Bunkyo-ku, Tokyo, 113-0033, Japan}
\affil[*]{Corresponding Author: lixiaoms@163.com}
\date{}
\maketitle

\abstract{In this study, simultaneous predictive distributions for independent Poisson observables were considered and the performance of predictive distributions was evaluated using the Kullback--Leibler (K--L) loss.
This study proposes a class of empirical Bayesian predictive distributions that dominate the Bayesian predictive distribution based on the Jeffreys prior. The K--L risk of the empirical Bayesian predictive distributions is demonstrated to be less than $1.04$ times the minimax lower bound.}

\keywords{
Predictive distribution; Kullback–-Leibler loss; Empirical Bayes; Minimaxity; Multivariate Poisson}

\section{Introduction}
The construction of accurate predictions is a fundamental problem in statistics. A reasonable approach is to construct a predictive distribution $q(y;x)$ to assign probabilities to possible future outcomes $y$ using the observed variables $x$. Therefore, the problem of constructing predictive distributions is highly important and has been studied in terms of various aspects \citep{aitchison1975goodness,komaki1996asymptotic,komaki2006,ghosh2019hierarchical}. As a representative discrete distribution, the Poisson distribution is commonly used to assume an integer data distribution. This study investigated the predictive distribution of Poisson observables.

The construction of the predictive distribution of Poisson observables is applicable to various fields. For example, different roads exist in a city, and the number of traffic accidents on each road per year is assumed to follow a Poisson distribution. The number of traffic accidents on each road in the following year can be predicted based on the number of traffic accidents in the past several years using the predictive distribution of Poisson observables. Prediction problems in various fields, such as sales and transportation, can also be formulated by constructing the predictive distribution of Poisson observables.

In the following, we assume that $x = (x_1,x_2,...,x_d)$ and $y = (y_1,y_2,...,y_d)$
are distributed according to the multivariate Poisson distributions,
$$p(x\mid\lambda)=\prod\limits_{i=1}^d p(x_i \mid \lambda_i)
=\exp\{-r(\lambda_1+\lambda_2+\cdots+\lambda_d)\}\frac{(r\lambda_1)^{x_1}}{x_1!}\cdots\frac{(r\lambda_d)^{x_d}}{x_d!}$$
and
$$p(y\mid\lambda)=\prod\limits_{i=1}^d p(y_i \mid \lambda_i)
=\exp\{-s(\lambda_1+\lambda_2+\cdots+\lambda_d)\}\frac{(s\lambda_1)^{y_1}}{y_1!}\cdots\frac{(s\lambda_d)^{y_d}}{y_d!},$$
respectively, where $r$ and $s$ are known positive real numbers. Let $\text{Po}(r\lambda)$ and $\text{Po}(s\lambda)$ denote the above Poisson distributions, respectively.
Here, $\lambda = (\lambda_1, \ldots, \lambda_d)$ is an unknown parameter.

We consider the problem of predicting the independent Poisson random variables
$y = (y_1,y_2,...,y_d)$ using the independent observations $x = (x_1,x_2,...,x_d)$. We adopt the Kullback--Leibler (K--L) loss of the predictive distribution $q(y;x)$, which is
$$D(p(y \mid \lambda),q(y;x))
=\sum\limits_yp(y\mid\lambda)
\log \frac{p(y\mid\lambda)}{q(y;x)}.$$ The K--L risk of the predictive distribution $q(y;x)$ on $\lambda$ is $\text{E}(D(p(y \mid \lambda),q(y;x))\mid \lambda).$

Numerous studies have been conducted on the estimation problem of the mean parameters of the multivariate Poisson distribution in the past century \citep{clevenson1975simultaneous,tsui1982simultaneous,ghosh1988simultaneous,chou1991simultaneous}. In contrast, studies on the predictive distribution problem of Poisson observables have recently emerged. \cite{komaki2004} proposed a class of shrinkage prior distributions $$\pi_{\alpha,\beta}(\lambda)\du\lambda_1\du\lambda_2\cdots \du\lambda_d\propto \frac{\lambda_1^{\beta_1-1}\lambda_2^{\beta_2-1}\cdots\lambda_d^{\beta_d-1}}
{(\lambda_1+\lambda_2+\cdots+\lambda_d)^\alpha}\du\lambda_1\du\lambda_2\cdots \du\lambda_d,$$ and the Bayesian predictive distribution based on $\pi_{\alpha=d/2-1,\beta=(1/2,\dots,1/2)}(\lambda)$ was shown to dominate that based on the Jeffreys prior. \cite{komaki2006class} proposed a class of proper priors and the Bayesian predictive distribution based on the proper priors was demonstrated to dominate that based on the Jeffreys prior.
More recently, \cite{hamura2020bayesian} studied the predictive distribution problem in a Poisson model with parametric restrictions. A class of asymptotic minimax Bayesian predictive distributions in sparse Poisson sequence models is presented in \cite{yano2021minimax}.

However, the construction of predictive distributions using the empirical Bayes approach has received little attention. A similar situation exists in predictive distribution studies of normal distributions. Although numerous studies have been conducted on Bayesian predictive distributions in normal models \citep{komaki2001shrinkage,brown2008,improve2011,matsuda2015singular}, relatively few works exist on empirical Bayesian predictive distributions. \cite{xu2011empirical} constructed a class of empirical Bayesian predictive distributions that were shown to dominate the Bayesian predictive distribution based on the Jeffreys prior, and were therefore minimax. Owing to the similarity between the Poisson and normal distributions in prediction theory \citep{komaki2006class}, we speculate that similar results can be obtained in the Poisson model, which is confirmed in this study. We use the empirical Bayes approach to construct a class of predictive distributions of Poisson observables, which are demonstrated to dominate the Bayesian predictive distribution based on the Jeffreys prior. Therefore, this study fills the gap in the research regarding the empirical Bayes prediction of Poisson observables.

In Section 2, we demonstrate that the Bayesian predictive distribution based on the Jeffreys prior is nearly minimax. More specifically, its K--L risk is less than $1.04$ times the minimax lower bound. In Section 3, we show that a class of empirical Bayesian predictive distributions dominates the Bayesian predictive distribution based on the Jeffreys prior. In Section 4, we compare the empirical Bayesian and Bayesian predictive distributions based on a shrinkage prior. Section 5 discusses different methods to design the value of the hyperparameter. The proofs of the main results are presented in Section 6.

\section{Bayesian predictive distribution under Jeffreys prior}

In this section, we consider the Bayesian predictive distribution based on the Jeffreys prior:
$$p_{\mathrm{J}}(y\mid x)=\frac{\int p(x,y\mid \lambda)\pi_{\mathrm{J}}(\lambda)\du \lambda}{\int p(x\mid \lambda)\pi_{\mathrm{J}}(\lambda)\du \lambda}=\frac{\int p(x\mid \lambda)p(y\mid \lambda)\pi_{\mathrm{J}}(\lambda)\du \lambda}{\int p(x\mid \lambda)\pi_{\mathrm{J}}(\lambda)\du \lambda},$$
where the Jeffreys prior $\pi_{\mathrm{J}}(\lambda)=\lambda_1^{-1/2}\lambda_2^{-1/2}\cdots\lambda_d^{-1/2}$. The analytical form of $p_{\mathrm{J}}(y\mid x)$ is presented in the following proposition.

\vspace{5pt}
\noindent
\begin{proposition}\label{Propo 1}
The Bayesian predictive distribution based on the Jeffreys prior is
$$p_{\mathrm{J}}(y\mid x)=\Big(\frac{r}{r+s}\Big)^{\sum_ix_i+d/2}\Big(\frac{s}{r+s}\Big)^{\sum_iy_i}\prod_{i=1}^d\frac{\Gamma(x_i+y_i+1/2)}{\Gamma(x_i+1/2)y_i!}.$$
\end{proposition}

First, we provide the upper bound for the maximum risk of $p_{\mathrm{J}}(y\mid x)$.
\vspace{5pt}
\noindent
\begin{theorem}\label{Theorem 1}
For any $\lambda$, the K--L risk of $p_{\mathrm{J}}(y\mid x)$ is less than $0.52d\log((r+s)/r)$.
\end{theorem}
\vspace{5pt}

Subsequently, we provide the lower bound for the minimax risk of predictive distributions.

\vspace{5pt}
\noindent
\begin{theorem}\label{Theorem 2}
For any predictive distribution $q(y;x)$ and positive number $\epsilon$, there exists $\lambda$ such that the K--L risk of $q(y;x)$ is greater than $0.5d\log((r+s)/r)-\epsilon$.
\end{theorem}
\vspace{5pt}

According to the two theorems, the upper bound of the K--L risk of $p_{\mathrm{J}}(y\mid x)$ is not greater than $1.04$ times the minimax lower bound. The minimax
risk divided by $0.5d\log((r + s)/r)$ lies in $[1, 1.04].$ The value $0.52$ in Theorem~\ref{Theorem 1} was obtained using a computer. We present the definition of a nearly minimax predictive distribution.

\vspace{5pt}
\noindent
\begin{definition}\label{Def 1}
A predictive distribution $q(y;x)$ is called nearly minimax if for any $\lambda$, the K--L risk of $q(y;x)$ is less than $1.04$ times the minimax lower bound.
\end{definition}
\vspace{5pt}

Hence, the Bayesian predictive distribution based on the Jeffreys prior is nearly minimax. Therefore, we are interested in the construction of a predictive distribution that is superior to $p_{\mathrm{J}}(y\mid x)$.

\section{A class of empirical Bayesian predictive distributions}

We describe the construction of the predictive distributions using the empirical Bayes approach. We consider an empirical Bayes model in which $x\sim\text{Po}(r\lambda)$, $y\sim\text{Po}(s\lambda)$, and $\lambda$ is distributed as a gamma prior:
\begin{equation}
    \lambda_i\sim \Gamma\Big(\frac{1}{2},\alpha\Big)= \lambda_i^{-1/2}\exp(-\lambda_i\alpha)\frac{\alpha^{1/2}}{\Gamma(1/2)},\ \text{iid}. \label{prior}
\end{equation}
The hyperparameter $\alpha$ is constructed using the observation $x$. Then, the empirical Bayesian predictive distribution under the gamma prior $\Gamma(\frac{1}{2},\alpha)$ is
$$\hat p_{\alpha}(y\mid x)=\frac{\int p(x\mid \lambda)p(y\mid \lambda)\prod_{i=1}^d\lambda_i^{-1/2}\exp(-\lambda_i\alpha)\du \lambda}{\int p(x\mid \lambda)\prod_{i=1}^d\lambda_i^{-1/2}\exp(-\lambda_i\alpha)\du \lambda}.$$
Although the form of the empirical Bayesian predictive distribution is the same as that of the Bayesian predictive distribution based on the gamma prior $\Gamma(\frac{1}{2},\alpha)$, in the empirical Bayesian predictive distribution $\hat p_{\alpha}(y\mid x)$, $\alpha$ changes according to the value of $x$, whereas in the Bayesian predictive distribution, $\alpha$ is a constant value. The analytical form of $\hat p_{\alpha}(y\mid x)$ is presented in the following proposition.

\vspace{5pt}
\noindent
\begin{proposition}\label{Propo 2}
The Bayesian predictive distribution based on the gamma prior \eqref{prior} is
$$\hat p_{\alpha}(y\mid x)=\Big(\frac{r+\alpha}{r+s+\alpha}\Big)^{\sum_ix_i+d/2}\Big(\frac{s}{r+s+\alpha}\Big)^{\sum_i y_i}\prod_{i=1}^d\frac{\Gamma(x_i+y_i+1/2)}{\Gamma(x_i+1/2)y_i!}.$$
\end{proposition}

Note that, under the empirical Bayes model, if $r$ is large and $d\ge3$,
\begin{align}
    &\text{E}\Big(\frac{r(d/2-1)}{\sum_{i=1}^d x_i+1}\:\Big|\;x\sim\text{Po}(r\lambda),\ \lambda_i\sim \Gamma\Big(\frac{1}{2},\alpha\Big)\Big) \notag\\
    &=\text{E}\Big(\frac{r(d/2-1)}{\sum_{i=1}^d r\lambda_i}\Big(1-\exp\Big(-\sum_{i=1}^d r\lambda_i\Big)\Big)\:\Big|\;\lambda_i\sim \Gamma\Big(\frac{1}{2},\alpha\Big)\Big) \notag\\
    &=\text{E}\Big(\frac{(d/2-1)}{\sum_{i=1}^d \lambda_i}\Big(1-\exp\Big(-\sum_{i=1}^d r\lambda_i\Big)\Big)\:\Big|\;\sum_{i=1}^d \lambda_i\sim \Gamma\Big(\frac{d}{2},\alpha\Big)\Big) \notag\\
    &\approx\text{E}\Big(\frac{d/2-1}{\sum_{i=1}^d \lambda_i}\:\Big|\;\sum_{i=1}^d \lambda_i\sim \Gamma\Big(\frac{d}{2},\alpha\Big)\Big) =\alpha. \notag
\end{align}
Therefore, a natural estimator of hyperparameter $\alpha$ is $r(d/2-1)/(\sum\nolimits_{i=1}^d x_i+1)$. We consider a general type of estimators $\hat\alpha=rb/(\sum\nolimits_{i=1}^d x_i+1)$, $0<b\le d-2$. We demonstrate that the empirical Bayesian predictive distribution $\hat p_{\hat\alpha}(y\mid x)$ dominates the Bayesian predictive distribution based on the Jeffreys prior.

\vspace{5pt}
\noindent
\begin{theorem}\label{Theorem 3}
If $d\ge3,$ $\alpha=rb/(\sum\nolimits_{i=1}^d x_i+1),$ and $0<b\le d-2$, $\hat p_{\alpha}(y\mid x)$ dominates $p_{\mathrm{J}}(y\mid x)$ and is thus nearly minimax. Furthermore, the risk difference between $\hat p_{\alpha}(y\mid x)$ and $p_{\mathrm{J}}(y\mid x)$ depends on $\lambda$ only through $\mu=\sum_{i=1}^d\lambda_i$.
\end{theorem}
\vspace{5pt}

\section{Comparison with Bayesian predictive distribution based on shrinkage prior}

In the previous section, we proposed a class of empirical Bayesian predictive distributions $\hat p_{\alpha}(y\mid x)$, where $\alpha=rb/(\sum\nolimits_{i=1}^d x_i+1)$ and $0<b\le d-2$. The empirical Bayesian predictive distributions $\hat p_{\alpha}(y\mid x)$ dominate the Bayesian predictive distribution $p_{\mathrm{J}}(y\mid x)$ based on the Jeffreys prior.

\begin{figure}[H]
	\centering
	\includegraphics[width=0.6\textwidth]{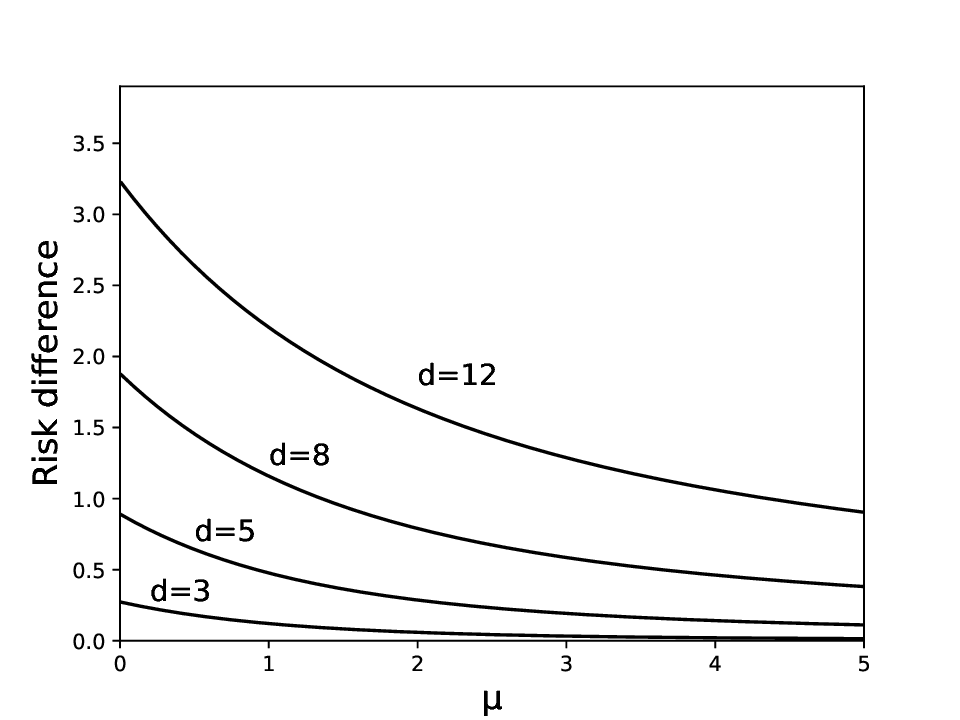}
	\caption{Risk difference between $p_{\mathrm{J}}(y\mid x)$ and $\hat p_{\alpha}(y\mid x)$ under different $\mu$ and $d$.}
        \label{fig:risk1}
\end{figure}

The K--L risk difference between predictive distributions $q_1$ and $q_2$ is defined as $R_\text{KL}(q_1)-R_\text{KL}(q_2),$ where $R_\text{KL}(q)$ denotes the K--L risk of $q$. Figure \ref{fig:risk1} shows the K--L risk differences between $p_{\mathrm{J}}(y\mid x)$ and $\hat p_{\alpha}(y\mid x)$ for the case $r=s=1$. Here, $\alpha=r(d/2-1)/(\sum\nolimits_{i=1}^d x_i+1).$ When $\mu$ is small, the risk difference is large. Therefore, the risk reduction that is offered by the empirical Bayesian predictive distribution is large if $\mu$ is small. Here, risk reduction offered by $q$ refers to the K--L risk difference between $p_{\mathrm{J}}(y\mid x)$ and $q$.

Next, $\hat p_{\alpha}(y\mid x)$ is compared with the Bayesian predictive distribution $p_{\mathrm{S}}(y\mid x)$ based on the shrinkage prior
$$ \pi_{\mathrm{S}}(\lambda)=(\lambda_1+\lambda_2+\cdots+\lambda_d)^{1-d/2}\lambda_1^{-1/2}\lambda_2^{-1/2}\cdots\lambda_d^{-1/2}.$$
We aim to compare the risk reductions that are offered by $p_{\mathrm{S}}(y\mid x)$ and $\hat p_{\alpha}(y\mid x)$.

\begin{figure}[H]
	\centering
	\includegraphics[width=1\textwidth]{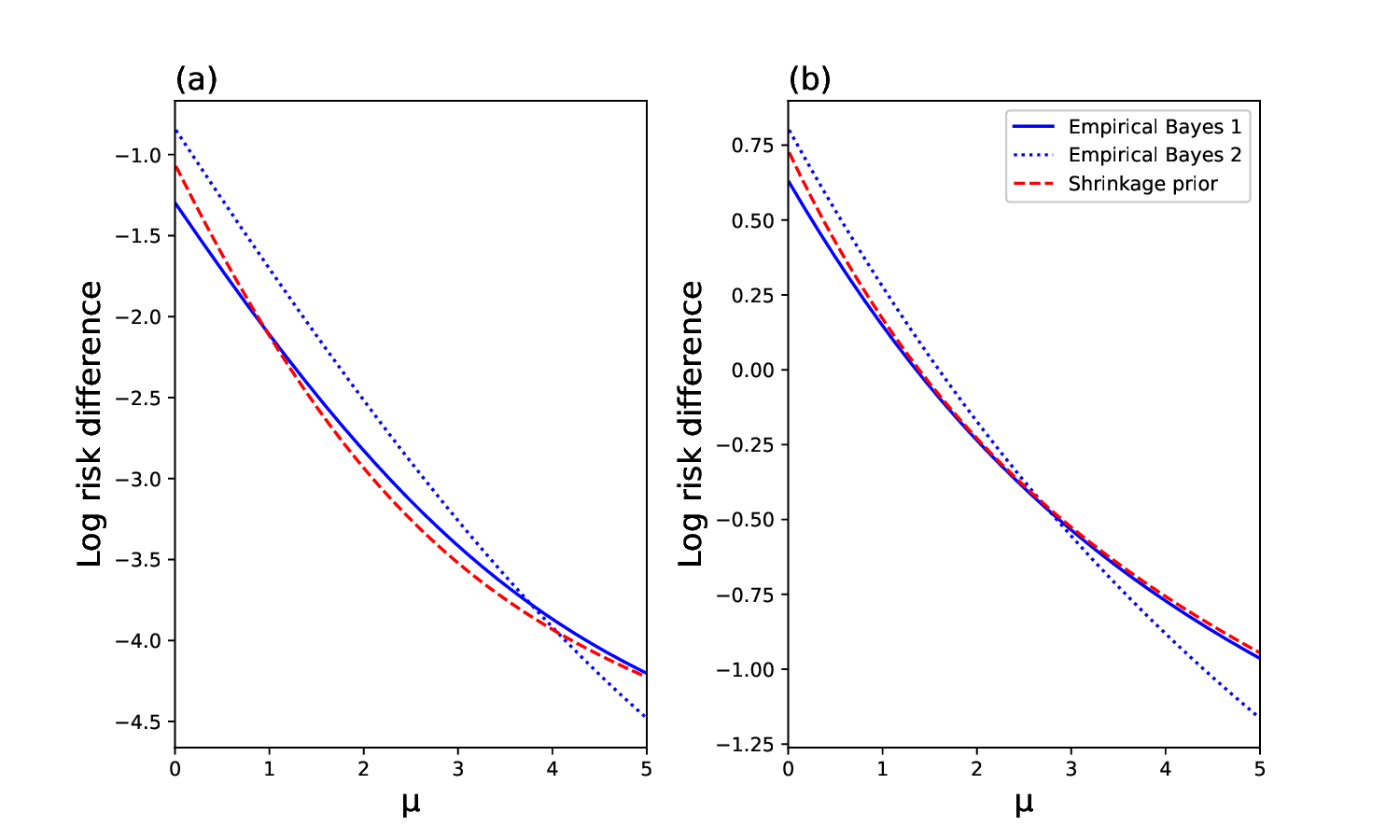}
	\caption{Log values of risk difference between $p_{\mathrm{J}}(y\mid x)$ and $\hat p_{\alpha}(y\mid x)$, and between $p_{\mathrm{J}}(y\mid x)$ and $p_{\mathrm{S}}(y\mid x)$ under different $\mu$ for (a) $d=3$ and (b) $d=8$.}
        \label{fig:risk2}
\end{figure}

We set $r=s=1$. Figure \ref{fig:risk2} shows the differences between the K--L risks of $p_{\mathrm{J}}(y\mid x)$ and empirical Bayesian predictive distributions $\hat p_{\alpha}(y\mid x)$, as well as between the K--L risks of $p_{\mathrm{J}}(y\mid x)$ and $p_{\mathrm{S}}(y\mid x)$. In the figure, empirical Bayes 1 denotes $\hat p_{\alpha_1}(y\mid x)$, where $\alpha_1=r(d/2-1)/(\sum\nolimits_{i=1}^d x_i+1)$, whereas empirical Bayes 2 denotes $\hat p_{\alpha_2}(y\mid x)$, where $\alpha_2=r(d-2)/(\sum\nolimits_{i=1}^d x_i+1)$. Subfigure (a) shows the results for the case $d=3$. It can be observed that when $\mu$ is smaller than $3$, the risk reduction offered by the empirical Bayesian predictive distribution $\hat p_{\alpha_2}(y\mid x)$ is the largest among the three predictive distributions. In contrast, when $\mu$ is larger than $4$, $\hat p_{\alpha_1}(y\mid x)$ and $p_{\mathrm{S}}(y\mid x)$ perform better than $\hat p_{\alpha_2}(y\mid x)$. $\hat p_{\alpha_1}(y\mid x)$ and $p_{\mathrm{S}}(y\mid x)$ perform similarly for each $\mu$. When $\mu$ is approximately $3$, $p_{\alpha_1}(y\mid x)$ outperforms $p_{\mathrm{S}}(y\mid x)$. Subfigure (b) shows the results for the case $d=8$, which are similar to those for $d=3$. $\hat p_{\alpha_2}(y\mid x)$ achieves the best performance for a small $\mu$ but worsens for a large $\mu$. $\hat p_{\alpha_1}(y\mid x)$ and $p_{\mathrm{S}}(y\mid x)$ perform similarly.

\section{Discussion}
This study proposes a class of empirical Bayesian predictive distributions of Poisson observables. The empirical Bayesian predictive distributions dominate the Bayesian predictive distribution based on the Jeffreys prior. Their K--L risk is demonstrated to be less than 1.04
times the minimax lower bound.

We used the approximate method of moments to determine the value of the hyperparameter $\alpha$. Here, the design of $\alpha$ is discussed from two other perspectives.
The first is maximum likelihood estimation (MLE). Under assumptions $x\sim\text{Po}(r\lambda)$ and $\lambda_i\sim\Gamma(\frac{1}{2},\alpha),\;\text{iid.}$, $$p(x\mid\alpha)=\prod_{i=1}^d\int\frac{(r\lambda_i)^{x_i}}{x_i!}e^{-r\lambda_i}\lambda_i^{-1/2}e^{-\alpha\lambda_i}\frac{\alpha^{1/2}}{\Gamma(1/2)}\du\lambda_i=\Big(\prod_{i=1}^d\frac{r^{x_i}\Gamma(x_i+1/2)}{x_i!\Gamma(1/2)}\Big)\alpha^{d/2}(r+\alpha)^{-\sum_ix_i-d/2}.$$ Maximizing $p(x\mid\alpha)$, the MLE $\hat\alpha=rd/(2\sum_{i=1}^dx_i)$ is obtained. 

The other is utilizing unbiased K--L risk estimate. \cite{george2021} proposed the unbiased estimate of the K--L risk of empirical predictive distributions in the normal model and designed the hyperparameters by minimizing the unbiased estimate. In the Poisson model of this study, using Proposition~\ref{Propo 2}, the K--L risk function of $\hat p_{\alpha}(y\mid x)$, which depends on $\alpha$ and $\lambda$, is
\begin{align}
&\text{E}\bigg(\log\Big(\prod_{i=1}^d\frac{(s\lambda_i)^{y_i}e^{-s\lambda_i}}{y_i!}\Big)-\log\Big(\Big(\frac{r+\alpha}{r+s+\alpha}\Big)^{\sum_ix_i+d/2}\Big(\frac{s}{r+s+\alpha}\Big)^{\sum_i y_i}\prod_{i=1}^d\frac{\Gamma(x_i+y_i+1/2)}{\Gamma(x_i+1/2)y_i!}\Big)\bigg)\notag\\
&=\sum_{i=1}^d\Big(s\lambda_i\log\lambda_i-s\lambda_i-\Big(r\lambda_i+\frac{1}{2}\Big)\log\Big(\frac{r+\alpha}{r+s+\alpha}\Big)+s\lambda_i\log(r+s+\alpha)-\text{E}\Big(\log\frac{\Gamma(x_i+y_i+\frac{1}{2})}{\Gamma(x_i+\frac{1}{2})}\Big)\Big). \label{KL-risk}
\end{align}
Similar to unbiased K--L risk estimate of estimators in Poisson model proposed by \cite{deledalle2017}, we ignore the terms in \eqref{KL-risk} that only depend on $\lambda$. Thus, we consider the remaining terms in \eqref{KL-risk}: $\sum_{i=1}^d\big(-(r\lambda_i+1/2)\log((r+\alpha)/(r+s+\alpha))+s\lambda_i\log(r+s+\alpha)\big).$ Therefore, $\alpha$ is chosen to minimize the unbiased estimate: $$U(\alpha)=\sum_{i=1}^d\Big(-\Big(x_i+\frac{1}{2}\Big)\log\Big(\frac{r+\alpha}{r+s+\alpha}\Big)+\frac{s}{r}x_i\log(r+s+\alpha)\Big).$$ $U(\alpha)$ achieves its minimum value at $\hat\alpha=rd/(2\sum_{i=1}^dx_i).$

Therefore, the choice of $\alpha$ obtained using the two methods is the same. However, $\hat\alpha=rd/(2\sum_{i=1}^dx_i)$ is not well-defined for the case of $\sum_{i=1}^dx_i=0$. Separately constructing a predictive distribution is needed for this case. Moreover, whether the corresponding empirical Bayesian predictive distribution dominates $p_{\mathrm{J}}(y\mid x)$ and whether it is nearly minimax require further study.

\section{Proofs}
\noindent
\textbf{Proof of Proposition~\ref{Propo 1}}.
\begin{align}
    p_{\mathrm{J}}(y\mid x)&=\frac{\int p(x\mid \lambda)p(y\mid \lambda)\pi_{\mathrm{J}}(\lambda)\du \lambda}{\int p(x\mid \lambda)\pi_{\mathrm{J}}(\lambda)\du \lambda} \notag\\
    &=\frac{\int \exp\{-(r+s)(\lambda_1+\lambda_2+\cdots+\lambda_d)\}\prod_{i=1}^d\frac{(r\lambda_i)^{x_i}}{x_i!}\frac{(s\lambda_i)^{y_i}}{y_i!} \lambda_i^{-1/2}\du \lambda}{\int \exp\{-r(\lambda_1+\lambda_2+\cdots+\lambda_d)\}\prod_{i=1}^d\frac{(r\lambda_i)^{x_i}}{x_i!} \lambda_i^{-1/2}\du \lambda} \notag\\
    &=\frac{\prod_{i=1}^d\int\exp(-(r+s)\lambda_i)\lambda_i^{x_i+y_i-1/2}\du\lambda_i}{\prod_{i=1}^d\int\exp(-r\lambda_i)\lambda_i^{x_i-1/2}\du\lambda_i}\prod_{i=1}^d\frac{s^{y_i}}{y_i!} \notag\\
    &=\Big(\frac{r}{r+s}\Big)^{\sum_ix_i+d/2}\Big(\frac{s}{r+s}\Big)^{\sum_iy_i}\prod_{i=1}^d\frac{\Gamma(x_i+y_i+1/2)}{\Gamma(x_i+1/2)y_i!}. \notag
\end{align}
\qed

\noindent
\textbf{Proof of Proposition~\ref{Propo 2}}.
\begin{align}
    \hat p_{\alpha}(y\mid x)&=\frac{\int p(x\mid \lambda)p(y\mid \lambda)\prod_{i=1}^d\lambda_i^{-1/2}\exp(-\lambda_i\alpha)\du \lambda}{\int p(x\mid \lambda)\prod_{i=1}^d\lambda_i^{-1/2}\exp(-\lambda_i\alpha)\du \lambda} \notag\\
    &=\frac{\int \exp\{-(r+s+\alpha)(\lambda_1+\lambda_2+\cdots+\lambda_d)\}\prod_{i=1}^d\frac{(r\lambda_i)^{x_i}}{x_i!}\frac{(s\lambda_i)^{y_i}}{y_i!} \lambda_i^{-1/2}\du \lambda}{\int \exp\{-(r+\alpha)(\lambda_1+\lambda_2+\cdots+\lambda_d)\}\prod_{i=1}^d\frac{(r\lambda_i)^{x_i}}{x_i!} \lambda_i^{-1/2}\du \lambda} \notag\\
    &=\frac{\prod_{i=1}^d\int\exp(-(r+s+\alpha)\lambda_i)\lambda_i^{x_i+y_i-1/2}\du\lambda_i}{\prod_{i=1}^d\int\exp(-(r+\alpha)\lambda_i)\lambda_i^{x_i-1/2}\du\lambda_i}\prod_{i=1}^d\frac{s^{y_i}}{y_i!} \notag\\
    &=\Big(\frac{r+\alpha}{r+s+\alpha}\Big)^{\sum_ix_i+d/2}\Big(\frac{s}{r+s+\alpha}\Big)^{\sum_iy_i}\prod_{i=1}^d\frac{\Gamma(x_i+y_i+1/2)}{\Gamma(x_i+1/2)y_i!}. \notag
\end{align}
\qed

\begin{figure}[H]
	\centering
	\includegraphics[width=.7\textwidth]{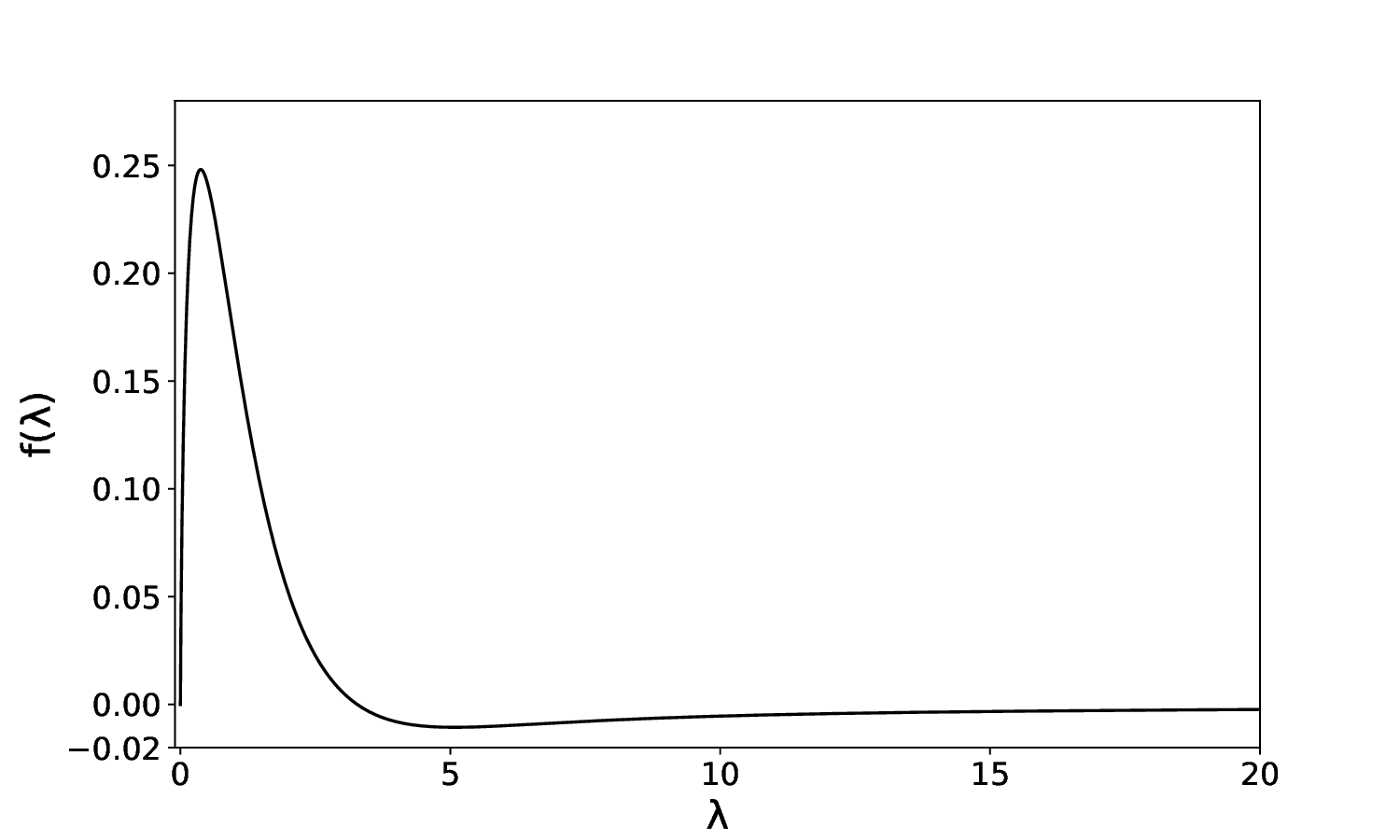}
	\caption{Graph of $f(\lambda)$ in $[0,20]$.}
        \label{fig:f}
\end{figure}

$f(\lambda)=\lambda\text{E}\big(\log((x+0.5)/\lambda)\mid x\sim\text{Po}(\lambda)\big)$ is defined. Figure \ref{fig:f} shows the graph of $f$ in the interval $(0,20]$. According to the numerical calculations, when $\lambda\in(0,20]$, $f$ achieves its minimum value around $5$, which is approximately $-0.011$.
The following lemmas are used for the proofs of the theorems. The proofs of the lemmas are presented in the Appendix.
\noindent
\begin{lemma}\label{Lemma 1}
For any $\lambda>0$, $f(\lambda)>-0.02$. Furthermore, $\lim_{\lambda\to\infty}f(\lambda)=0$.
\end{lemma}

\noindent
\begin{lemma}\label{Lemma 2}
For any $x>0$, $t>0$, and $s>0$, $$-(x+t+1)\log\Big(1+\frac{s}{1+s}\frac{-2t}{x+2t+1}\Big)-sx\log\Big(1+\frac{1}{1+s}\frac{2t}{x}\Big)>0.$$
\end{lemma}

\noindent
\textbf{Proof of Theorem~\ref{Theorem 1}}.

According to Proposition \ref{Propo 1}, the K--L risk $\text{E}\big(D(p(y \mid \lambda),p_{\mathrm{J}}(y\mid x))\big)$ is given by
\begin{align}
    &\text{E}\Big(\log p(y\mid \lambda)-\log p_{\mathrm{J}}(y\mid x)\:\Big|\;x\sim \text{Po}(r\lambda),\ y\sim \text{Po}(s\lambda)\Big)\notag\\
    &=\text{E}\Big(\log p(y\mid \lambda)-\log\Big(\Big(\frac{r}{r+s}\Big)^{\sum_ix_i+d/2}\Big(\frac{s}{r+s}\Big)^{\sum_iy_i}\prod_{i=1}^d\frac{\Gamma(x_i+y_i+1/2)}{\Gamma(x_i+1/2)y_i!}\Big)\:\Big|\;\lambda \Big)\notag\\
    &=\text{E}\Big(\sum_{i=1}^d\big(y_i\log(s\lambda_i)-s\lambda_i\big)-\Big(\sum_ix_i+\frac{d}{2}\Big)\log\Big(\frac{r}{r+s}\Big)-\Big(\sum_iy_i\Big)\log\Big(\frac{s}{r+s}\Big)-\log\prod_{i=1}^d\frac{\Gamma(x_i+y_i+\frac{1}{2})}{\Gamma(x_i+\frac{1}{2})}\Big)\notag\\
    &=\sum_{i=1}^d\biggl(-s\lambda_i+s\lambda_i\log(s\lambda_i)-\Big(r\lambda_i+\frac{1}{2}\Big)\log\Big(\frac{r}{r+s}\Big)-s\lambda_i\log\Big(\frac{s}{r+s}\Big)\notag\\
    &\qquad-\text{E}\Bigl(\log\Gamma\Big(x_i+y_i+\frac{1}{2}\Big)-\log\Gamma\Big(x_i+\frac{1}{2}\Big)\:\Big|\;\lambda \Bigr) \biggr) \label{appen-1}
\end{align}
Considering the function
$$F(t):=\sum_{i=1}^d\biggl(t\lambda_i\log\lambda_i+\frac{1}{2}\log t+\lambda_i(t\log t-t)-\text{E}\Bigl(\log\Gamma\Big(x+\frac{1}{2}\Big)\:\Big|\;x\sim \text{Po}(t\lambda_i)\Bigr)\biggr),$$
because $x_i+y_i\sim\text{Po}((r+s)\lambda_i),$ the K--L risk \eqref{appen-1} is equal to $F(r+s)-F(r)=\int_{r}^{r+s} F^{\prime}(t)\du t.$
We have
\begin{align}
F^{\prime}(t)&=\sum_{i=1}^d\biggl(\lambda_i\log\lambda_i+\frac{1}{2t}+\lambda_i\log t-\Bigl(\sum_{x\ge0}\log\Gamma\Big(x+\frac{1}{2}\Big)\frac{(t\lambda_i)^x}{x!}\exp(-t\lambda_i)\Bigr)^{\prime}\biggr) \notag\\
&=\sum_{i=1}^d\biggl(\frac{1}{2t}+\lambda_i\biggl(\log(t\lambda_i)-\sum_{x\ge1}\log\Gamma\Big(x+\frac{1}{2}\Big)\frac{(t\lambda_i)^{x-1}}{(x-1)!}e^{-t\lambda_i}+\sum_{x\ge0}\log\Gamma\Big(x+\frac{1}{2}\Big)\frac{(t\lambda_i)^x}{x!}e^{-t\lambda_i}\biggr)\biggr) \notag\\
&=\sum_{i=1}^d\biggl(\frac{1}{2t}+\lambda_i\biggl(\log(t\lambda_i)-\sum_{x\ge0}\Big(\log\frac{\Gamma(x+1.5)}{\Gamma(x+0.5)}\Big)\frac{(t\lambda_i)^x}{x!}e^{-t\lambda_i}\biggr)\biggr) \notag\\
&=\sum_{i=1}^d\biggl(\frac{1}{2t}-\text{E}\Bigl(\lambda_i\log\Big(\frac{x+0.5}{t\lambda_i}\Big)\:\Big|\;x\sim \text{Po}(t\lambda_i)\Bigr) \biggr)= \sum_{i=1}^d\Big(\frac{1}{2t}-\frac{1}{t}f(t\lambda_i)\Big).  \label{diff-F}
\end{align}
Therefore, the K--L risk \eqref{appen-1} is equal to \begin{align}
\int_{r}^{r+s} \sum_{i=1}^d\Big(\frac{1}{2t}-\frac{1}{t}f(t\lambda_i)\Big)\du t.    \label{appen-1.1}
\end{align}
From Lemma \ref{Lemma 1}, $f(t\lambda_i)>-0.02.$ Thus, the K--L risk \eqref{appen-1} is less than $0.52d\log((r+s)/r)$.

\qed

\noindent
\textbf{Proof of Theorem~\ref{Theorem 2}}.

We only need to show that $0.5d\log((r+s)/r)$ is the Bayes risk limit of a sequence of Bayes rules $p_{\pi_n}$ with $$\pi_n(\lambda)=\prod_{i=1}^d\lambda_i^{-1/2}\exp(-\frac{\lambda_i}{n})\frac{1}{n^{1/2}\Gamma(1/2)}.$$
The Bayes risk of $p_{\pi_n}$ is equal to
\begin{align}
\text{E}\Big(\text{E}\Big(\log \frac{p(y\mid\lambda)}{p_{\mathrm{J}}(y\mid x)}\Big) \:\Big|\;\lambda\sim\pi_n\Big)+\text{E}\Big(\text{E}\Big(\log \frac{p_{\mathrm{J}}(y\mid x)}{p_{\pi_n}(y\mid x)}\Big)\:\Big|\;\lambda\sim\pi_n\Big). \label{appen-2}
\end{align}

We first show that the left term in \eqref{appen-2} converges to $0.5d\log((r+s)/r)$ when $n\to\infty$.
Using \eqref{appen-1.1} in the proof of Theorem \ref{Theorem 1}, the left term in \eqref{appen-2} is equal to
\begin{align}
   &\text{E}\bigg(\int_{r}^{r+s} \sum_{i=1}^d\Big(\frac{1}{2t}-\frac{1}{t}f(t\lambda_i) \Bigr) \du t\:\Big|\;\lambda\sim\pi_n\bigg) = 0.5d\log\Big(\frac{r+s}{r}\Big)-\sum_{i=1}^d\int_{r}^{r+s}\frac{1}{t}\text{E}\bigg(f(t\lambda_i)\:\Big|\;\lambda_i\sim\Gamma\Big(\frac{1}{2},\frac{1}{n}\Big)\bigg)\du t \label{appen-3}.
\end{align}
According to $\lim_{\lambda\to\infty}f(\lambda)=0$ from Lemma \ref{Lemma 1}, \eqref{appen-3} converges to $0.5d\log((r+s)/r)$ when $n\to\infty$.

We then show that the right term in \eqref{appen-2} converges to 0 when $n\to\infty$. From Proposition \ref{Propo 2}, we obtain
$$p_{\pi_n}(y\mid x)=\Big(\frac{r+1/n}{r+s+1/n}\Big)^{\sum_ix_i+d/2}\Big(\frac{s}{r+s+1/n}\Big)^{\sum_i y_i}\prod_{i=1}^d\frac{\Gamma(x_i+y_i+1/2)}{\Gamma(x_i+1/2)y_i!}.$$
Therefore, $$\frac{p_{\mathrm{J}}(y\mid x)}{p_{\pi_n}(y\mid x)}=\Big(\frac{r}{r+1/n}\Big)^{\sum_ix_i+d/2}\Big(\frac{r+s+1/n}{r+s}\Big)^{\sum_i(x_i+y_i)+d/2}.$$ When $n\to\infty$, the right term in \eqref{appen-2} is equal to
\begin{align}
&\text{E}\bigg(\text{E}\Big(\Big(\sum_ix_i+d/2\Big)\log\Big(\frac{r}{r+1/n}\Big)+\Big(\sum_i(x_i+y_i)+d/2\Big)\log\Big(\frac{r+s+1/n}{r+s}\Big)\Big) \:\Big|\;\lambda\sim\pi_n\bigg) \notag\\
   &=\text{E}\Big(\big(r\mu+d/2\big)\log\Big(\frac{r}{r+1/n}\Big)+\big((r+s)\mu+d/2\big)\log\Big(\frac{r+s+1/n}{r+s}\Big)\:\Big|\;\mu\sim\Gamma\Big(\frac{d}{2},\frac{1}{n}\Big)\Big)  \notag\\
   &\to\text{E}\Big(r\mu\log\Big(\frac{r}{r+1/n}\Big)+(r+s)\mu\log\Big(\frac{r+s+1/n}{r+s}\Big)\:\Big|\;\mu\sim\Gamma\Big(\frac{d}{2},\frac{1}{n}\Big)\Big)  \notag\\
   &= dn/2\Big(r\log\Big(\frac{r}{r+1/n}\Big)+(r+s)\log\Big(\frac{r+s+1/n}{r+s}\Big)\Big) \to 0, \notag
\end{align}
where $\mu=\sum_{i=1}^d\lambda_i.$

Therefore, \eqref{appen-2} converges to $0.5d\log((r+s)/r)$ when $n\to\infty$, which completes the proof.

\qed

\noindent
\textbf{Proof of Theorem~\ref{Theorem 3}}.

From Propositions \ref{Propo 1} and \ref{Propo 2}, when $\alpha=rb/(\sum\nolimits_{i=1}^d x_i+1)$, the K--L risk difference between $p_{\mathrm{J}}(y\mid x)$ and $\hat p_{\alpha}(y\mid x)$ is given by
\begin{align}
&\text{E}\Big(
\log \frac{\hat p_{\alpha}(y\mid x)}{p_{\mathrm{J}}(y\mid x)}\:\Big|\;x\sim\text{Po}(r\lambda),\ y\sim\text{Po}(s\lambda)\Big) \notag\\
&=\text{E}\bigg(\Big(\sum_ix_i+d/2\Big)\log\Big(\frac{r+\alpha}{r}\Big)-\Big(\sum_i(x_i+y_i)+d/2\Big)\log\Big(\frac{r+s+\alpha}{r+s}\Big)\:\Big|\;x\sim\text{Po}(r\lambda),\ y\sim\text{Po}(s\lambda)\bigg) \notag\\
&=\text{E}\biggl(\Big(\sum_ix_i+d/2\Big)\log\Big(\frac{\sum_ix_i+1+b}{\sum_ix_i+1}\Big)-\Big(\sum_ix_i+d/2\Big)\log\Big(\frac{s+r(\sum_ix_i+1+b)/(\sum_ix_i+1)}{s+r}\Big) \notag\\
&\quad-\Big(\sum_iy_i\Big)\log\Big(\frac{s+r(\sum_ix_i+1+b)/(\sum_ix_i+1)}{s+r}\Big)
\:\Big|\;x\sim \text{Po}(r\lambda),\ y\sim \text{Po}(s\lambda)\biggr) \notag\\
&=\text{E}\biggl(\Big(X+\frac{d}{2}\Big)\log\Big(\frac{(X+1+b)/(X+1)}{(s+r\frac{X+1+b}{X+1})/(s+r)}\Big) -Y\log\Big(\frac{s+\frac{r(X+1+b)}{X+1}}{s+r}\Big)
\:\Big|\;X\sim \text{Po}\Big(r\sum_i\lambda_i\Big),\ Y\sim \text{Po}\Big(s\sum_i\lambda_i\Big)\biggr), \label{appen-4}
\end{align}
where $X=\sum_ix_i,$ $Y=\sum_iy_i.$
Note that for any function $g(X)$, 
\begin{align*}
&\text{E}\Big(Yg(X)\:\Big|\;X\sim \text{Po}\Big(r\sum_i\lambda_i\Big),\ Y\sim \text{Po}\Big(s\sum_i\lambda_i\Big)\Big)\\
&=s\Big(\sum_i\lambda_i\Big)\sum_{X\ge0}g(X)\frac{(r\sum_i\lambda_i)^X}{X!}e^{-r\sum_i\lambda_i}=\frac{s}{r}\sum_{X\ge0}g(X)(X+1)\frac{(r\sum_i\lambda_i)^{X+1}}{(X+1)!}e^{-r\sum_i\lambda_i}\\
&=\frac{s}{r}\text{E}\Big(Xg(X-1)\:\Big|\;X\sim \text{Po}\Big(r\sum_i\lambda_i\Big)\Big).
\end{align*}
Thus, the K--L risk difference \eqref{appen-4} is equal to
\begin{align}
&\text{E}\biggl(-\Big(X+\frac{d}{2}\Big)\log\Big(\frac{(s+r\frac{X+1+b}{X+1})/(s+r)}{(X+1+b)/(X+1)}\Big)
-\frac{s}{r}X\log\Big(\frac{s+\frac{r(X+b)}{X}}{s+r}\Big)
\:\Big|\;X\sim \text{Po}\Big(r\sum_i\lambda_i\Big)\biggr) \notag\\
&=\text{E}\biggl(-\Big(X+\frac{d}{2}\Big)\log\Big(1+\frac{s}{r+s}\frac{-b}{X+1+b}\Big)-\frac{s}{r}X\log\Big(1+\frac{r}{r+s}\frac{b}{X}\Big)\:\Big|\;X\sim \text{Po}\Big(r\sum_i\lambda_i\Big)\biggr), \label{appen-5}
\end{align}
which depends on $\lambda$ only through $\mu=\sum_{i=1}^d\lambda_i$. According to Lemma \ref{Lemma 2}, we obtain
$$-\Big(x+\frac{b}{2}+1\Big)\log\Big(1+\frac{s/r}{1+s/r}\frac{-b}{x+b+1}\Big)-\frac{s}{r}x\log\Big(1+\frac{1}{1+s/r}\frac{b}{x}\Big)>0.$$
In combination with $d/2\ge b/2+1$, we obtain
$$-\Big(X+\frac{d}{2}\Big)\log\Big(1+\frac{s}{r+s}\frac{-b}{X+1+b}\Big)-\frac{s}{r}X\log\Big(1+\frac{r}{r+s}\frac{b}{X}\Big)>0,\ \forall X\ge0.$$
Thus, \eqref{appen-5} is positive. Therefore, \eqref{appen-4} is positive, which completes the proof.

\qed

\section*{Acknowledgments}
I am grateful for the support from the China Scholarship Council. I thank Yangkendi Deng
and Haokun Li for their helpful comments on the proofs. I also thank Fumiyasu Komaki and Takeru Matsuda for their helpful comments. I am grateful to the anonymous referees for their constructive comments.

\begin{appendices}
\section{Proof of Lemma~\ref{Lemma 1}.}
\noindent
\textit{\textbf{Proof of part 1.}}
First, we prove that $f(\lambda)=\lambda\text{E}\big(\log((x+0.5)/\lambda)\mid x\sim\text{Po}(\lambda)\big)>-0.02,\ \forall\lambda>0$ in two cases: $\lambda\le1$ and $\lambda>1$. We present the outline of the proof's flow as follows:

When $\lambda\le1$, we prove $f(\lambda)>0$ using $\text{E}\big(\log((x+0.5)/\lambda)\mid x\sim\text{Po}(\lambda)\big)>\log(0.5/\lambda)P(x=0)+\log(1.5/\lambda)P(x\ge1).$

When $\lambda>1$, we define the derivative of $f(\lambda)/\lambda$ as $g(\lambda).$ We derive a lower bound \eqref{low-2} and an upper bound \eqref{upp-2} for $g(\lambda).$ We used a computer to verify that $f(3)>0$, $f(4)>-0.0082$, and $f(5)>-0.011.$ Using these values and upper and lower bounds for $g(\lambda)$, we can obtain $f(\lambda)>-0.02.$ 

The details of each case are presented below.

\vspace{0.2cm}
\textbf{Case 1:} $\lambda\le1$.

When $\lambda\le1/2$, $(x+0.5)/\lambda\ge1$. Thus, $f(\lambda)\ge0$. 

When $1/2<\lambda<1$, $\text{E}\big(\log((x+0.5)/\lambda)\mid x\sim\text{Po}(\lambda)\big)>\log(0.5/\lambda)P(x=0)+\log(1.5/\lambda)P(x\ge1)=\log(0.5/\lambda)e^{-\lambda}+\log(1.5/\lambda)(1-e^{-\lambda})=\log(1.5/\lambda)-(\log3)e^{-\lambda},$ which is positive because $(\log(1.5/\lambda)-(\log3)e^{-\lambda})'=-1/\lambda+(\log3)e^{-\lambda}<-1+(\log3)e^{-1/2}<0$ and $\log(1.5)-(\log3)e^{-1}>0.$

\vspace{0.2cm}
\textbf{Case 2:} $\lambda>1$.

Let $g(\lambda)$ denote the derivative of $f(\lambda)/\lambda$.
\begin{align}
g(\lambda)&=\Big(\sum_{x=0}^{\infty}\log\Big(\frac{x+0.5}{\lambda}\Big)e^{-\lambda}\frac{\lambda^x}{x!}\Big)'\notag\\
&=\sum_{x=0}^{\infty}\Big(\log\Big(\frac{x+0.5}{\lambda}\Big)e^{-\lambda}\frac{\lambda^{x-1}x}{x!}-\log\Big(\frac{x+0.5}{\lambda}\Big)e^{-\lambda}\frac{\lambda^{x}}{x!}-e^{-\lambda}\frac{\lambda^{x-1}}{x!}\Big) \notag\\
&=\sum_{x=0}^{\infty}\Big(\log\Big(\frac{x+1.5}{x+0.5}\Big)e^{-\lambda}\frac{\lambda^x}{x!}\Big)-1/\lambda=\text{E}\Big(\log\Big(\frac{x+1.5}{x+0.5}\Big)\Big)-1/\lambda.\label{diff-g}
\end{align}
For any $t\ge0$, note that the Taylor’s formula
\begin{align}
\log\Big(\frac{t+1.5}{t+0.5}\Big)=\log\Big(1+\frac{1}{2(t+1)}\Big)-\log\Big(1-\frac{1}{2(t+1)}\Big)=\sum_{k=1}^\infty\frac{2}{2k-1}\Big(\frac{1}{2(t+1)}\Big)^{2k-1}.  \label{taylor}
\end{align}
Thus, $\log(t+1.5)-\log(t+0.5)>1/(t+1).$ Using \eqref{diff-g}, we obtain
\begin{align}
g(\lambda)&=\text{E}\Big(\log\Big(\frac{x+1.5}{x+0.5}\Big)\Big)-1/\lambda=(\log3)P(x=0)+\text{E}\Big(\log\Big(\frac{x+1.5}{x+0.5}\Big)\mathbf{1}(x\ge1)\Big)-1/\lambda\notag\\
&>1.09P(x=0)+\text{E}\Big(\frac{1}{x+1}\mathbf{1}(x\ge1)\Big)-1/\lambda=0.09P(x=0)+\text{E}\Big(\frac{1}{x+1}\Big)-1/\lambda. \label{low-1}
\end{align}
Because $\text{E}(\lambda/(x+1))=1-e^{-\lambda},$ from \eqref{low-1},
we obtain a lower bound of $g(\lambda)$:
\begin{align}
    g(\lambda)>0.09e^{-\lambda}-e^{-\lambda}\lambda^{-1}. \label{low-2}
\end{align}
Using the Taylor’s formula \eqref{taylor}, for any $t\ge2$,
\begin{align}
\log\Big(\frac{t+1.5}{t+0.5}\Big)<\frac{1}{t+1}+\sum_{k=2}^\infty\frac{2^{2-2k}}{2k-1}(t+1)^{-3}=\frac{1}{t+1}+\frac{\log3-1}{(t+1)^3}<\frac{1}{t+1}+\frac{0.26}{(t+1)(t+2)(t+3)}. \label{upp-1}
\end{align}
Because $\log(2.5/1.5)<0.5+0.26/24$, combining \eqref{diff-g} and \eqref{upp-1}, we obtain
\begin{align}
    g(\lambda)&=\text{E}\Big(\log\Big(\frac{x+1.5}{x+0.5}\Big)\Big)-\frac{1}{\lambda}\notag=(\log3)P(x=0)+\text{E}\Big(\log\Big(\frac{x+1.5}{x+0.5}\Big)\mathbf{1}(x\ge1)\Big)-\frac{1}{\lambda}\notag\\
    &<(\log3-1-0.26/6)P(x=0)+\text{E}\Big(\frac{1}{x+1}+\frac{0.26}{(x+1)(x+2)(x+3)}\Big)-\frac{1}{\lambda}\notag.
\end{align}
Because $\text{E}(\lambda/(x+1))=1-e^{-\lambda}$ and $\text{E}((x+1)^{-1}(x+2)^{-1}(x+3)^{-1})<\lambda^{-3}$, we get an upper bound of $g(\lambda)$:
\begin{align}
    g(\lambda)<0.06e^{-\lambda}-e^{-\lambda}\lambda^{-1}+0.26\lambda^{-3}. \label{upp-2}
\end{align}

Using a computer, we can calculate the value of function $$L(\lambda)=\sum_{x=0}^{20}\log(x+0.5)\frac{\lambda^x}{x!}\exp(-\lambda)\lambda-\lambda\log\lambda $$ for $\lambda=3,4,5.$ We only calculate $x\le20$ to calculate only a finite number of terms. The code for the calculation and the analysis of potential numerical errors are available at \href{https://github.com/lixiaoms/EB-Poisson}{https://github.com/lixiaoms/EB-Poisson}. We obtained $f(3)>L(3)>0$, $f(4)>L(4)>-0.0082$, and $f(5)>L(5)>-0.011$. Next, we use these inequalities and the upper and lower bounds of $g(\lambda)$ to prove that $f(\lambda)>-0.02.$ We prove it in five cases as follows. The selection of $3$, $4$, $5$, and $7$ as the boundaries for different cases is because the inequality discussed in each case holds in the corresponding interval, and the lower bounds of $f(3)$, $f(4)$, and $f(5)$ are used.

(1) Case of $\lambda\ge7$. From \eqref{upp-2}, $g(t)<0.06e^{-t}+0.26t^{-3}.$
Because $g(\lambda)=(f(\lambda)/\lambda)'$ and $\lim_{\lambda\to\infty}f(\lambda)=0$ (the proof is presented in the second part of Appendix A), we have
$$f(\lambda)/\lambda=-\int_{\lambda}^{\infty}g(t)\du t>-\int_{\lambda}^{\infty}(0.06e^{-t}+0.26t^{-3})\du t=-0.06e^{-\lambda}-0.13\lambda^{-2}.$$
Thus, $f(\lambda)>-0.06e^{-\lambda}\lambda-0.13/\lambda\ge-0.06\times e^{-7}\times7-0.13/7>-0.02.$

(2) Case of $\lambda\in[5,7].$ From \eqref{low-2}, when $t>5$, $g(t)>0.09e^{-t}-e^{-t}t^{-1}>-0.11e^{-t}.$ Thus,
$$f(\lambda)/\lambda=f(5)/5+\int_{5}^{\lambda}g(t)\du t>-0.011/5-\int_{5}^{\lambda}0.11e^{-t}\du t>-0.00295+0.11e^{-\lambda}.$$
Thus, $f(\lambda)>-0.00295\lambda+0.11e^{-\lambda}\lambda\ge-0.00295\times7+0.11e^{-7}\times7>-0.02.$

(3) Case of $\lambda\in[4,5].$ From \eqref{low-2}, when $t>4$, $g(t)>0.09e^{-t}-e^{-t}t^{-1}>-0.16e^{-t}.$ Thus, 
$$f(\lambda)/\lambda=f(4)/4+\int_{4}^{\lambda}g(t)\du t>-0.0082/4-\int_{4}^{\lambda}0.16e^{-t}\du t>-0.005+0.16e^{-\lambda}.$$
Thus, $f(\lambda)>-0.005\lambda+0.16e^{-\lambda}\lambda\ge-0.005\times5+0.16e^{-5}\times5>-0.02.$

(4) Case of $\lambda\in[3,4].$ From \eqref{upp-2}, when $t<4$, $g(t)<0.06e^{-t}-e^{-t}t^{-1}+0.26t^{-3}<-0.19e^{-t}+0.26t^{-3}.$ When $t\in(3,4)$, using $e^{-t}t^3>e^{-4}4^3>1$, we have $g(t)<-0.19t^{-3}+0.26t^{-3}=0.07t^{-3}$. Thus,
$$f(\lambda)/\lambda=f(4)/4-\int_{\lambda}^{4}g(t)\du t>-0.0082/4-\int_{\lambda}^{4}0.07t^{-3}\du t>-0.035\lambda^{-2}.$$
Thus, $f(\lambda)>-0.035\lambda^{-1}>-0.02.$

(5) Case of $\lambda\in[1,3].$ From \eqref{upp-2}, when $t\in(1,3)$, $g(t)<0.06e^{-t}-e^{-t}t^{-1}+0.26t^{-3}<0.2e^{-t}t^{-1}-e^{-t}t^{-1}+0.26t^{-3}=-0.8e^{-t}t^{-1}+0.26t^{-3}.$ Because $e^{-t}t^2>\max(e^{-3}\times3^2,e^{-1})>0.36$ when $t\in(1,3)$, we obtain $g(t)<-0.8\times0.36t^{-3}+0.26t^{-3}<0.$ Therefore, $f(\lambda)/\lambda$ is decreasing in $[1,3].$ Using $f(3)>0$, we obtain $f(\lambda)>0$ for any $\lambda\in[1,3].$

\vspace{0.2cm}
\noindent
\textit{\textbf{Proof of part 2.}} Subsequently, we prove that $\lim_{\lambda\to\infty}f(\lambda)=0$.

First, we prove $\lim\inf_{\lambda\to\infty}f(\lambda)\ge0$. For any given $\epsilon>0$, there exists $\delta\in(0,\ 0.1)$ such that $\log(1+t)\ge t-(0.5+\epsilon)t^2,\ \forall t\ge-2\delta.$ Without loss of generality, we assume $\lambda>1/\delta$. Therefore, by setting $t=(x+0.5-\lambda)/\lambda,$ we obtain
\begin{align}
f(\lambda)&=\lambda\text{E}\Big(\log\Big(\frac{x+0.5}{\lambda}\Big)\mathbf{1}(x<(1-\delta)\lambda)\Big)+\lambda\text{E}\Big(\log\Big(\frac{x+0.5}{\lambda}\Big)\mathbf{1}(x\ge(1-\delta)\lambda)\Big)\notag\\
&\ge\lambda\text{E}\Big(\log\Big(\frac{x+0.5}{\lambda}\Big)\mathbf{1}(x<(1-\delta)\lambda)\Big)+\lambda\text{E}\Big(\Big(\frac{x+0.5-\lambda}{\lambda}-(0.5+\epsilon)\Big(\frac{x+0.5-\lambda}{\lambda}\Big)^2\Big)\mathbf{1}(x\ge(1-\delta)\lambda)\Big)\notag\\
&\ge\lambda\text{E}\Big(\log\Big(\frac{x+0.5}{\lambda}\Big)\mathbf{1}(x<(1-\delta)\lambda)\Big)+\lambda\text{E}\Big((x+0.5-\lambda)/\lambda-(0.5+\epsilon)(x+0.5-\lambda)^2/\lambda^2\Big).\label{lem1-2.1}
\end{align}
Using Chernoff bound for Poisson distribution, we obtain $$P(x<(1-\delta)\lambda)\le\frac{(e\lambda)^{(1-\delta)\lambda}e^{-\lambda}}{((1-\delta)\lambda)^{(1-\delta)\lambda}}.$$ When $x<(1-\delta)\lambda$, $|\log((x+0.5)/\lambda)|\le\log(2\lambda).$ Thus, the logarithm of the absolute value of the first term in the \eqref{lem1-2.1} is not greater than
\begin{align*}
&\log\Big(\lambda\log(2\lambda)P(x<(1-\delta)\lambda)\Big)\le \log\Big(\lambda\log(2\lambda)\frac{(e\lambda)^{(1-\delta)\lambda}e^{-\lambda}}{((1-\delta)\lambda)^{(1-\delta)\lambda}}\Big)\\
&=\log(\lambda\log(2\lambda))+(1-\delta)\lambda\log\lambda+(1-\delta)\lambda-\lambda-(1-\delta)\lambda\log\big((1-\delta)\lambda\big)\\
&=\big(-(1-\delta)\log(1-\delta)-\delta\big)\lambda+\mathrm{o}(\lambda)\to -\infty
\end{align*}
when $\lambda\to\infty.$ Thus, the first term of \eqref{lem1-2.1} converges to $0$ when $\lambda\to\infty.$ Because $\text{E}((x+0.5-\lambda)^2)=\lambda+0.25$, the second term of \eqref{lem1-2.1} is $-\epsilon-(0.5+\epsilon)0.25/\lambda$. Thus, $\eqref{lem1-2.1}\to-\epsilon$ when $\lambda\to\infty.$ Thus, $\lim\inf_{\lambda\to\infty}f(\lambda)\ge-\epsilon$. Because $\epsilon$ is an arbitrary positive value, $\lim\inf_{\lambda\to\infty}f(\lambda)\ge0.$

Next, we prove $\lim\sup_{\lambda\to\infty}f(\lambda)\le0.$ Note that $\log(1+t)\le t-t^2/2+t^3/3,\ \forall t.$ Thus, using $\text{E}((x-\lambda)^3)=\lambda$,
\begin{align}
f(\lambda)&=\lambda\text{E}\Big(\log\Big(\frac{x+0.5}{\lambda}\Big)\;\Big|\;x\sim\text{Po}(\lambda)\Big)\notag\\
&\le \lambda\text{E}\Big(\frac{x+0.5-\lambda}{\lambda}-\frac{(x+0.5-\lambda)^2}{2\lambda^2}+\frac{(x+0.5-\lambda)^3}{3\lambda^3}\;\Big|\;x\sim\text{Po}(\lambda)\Big)\notag\\
&=0.5-(\lambda+0.5^2)/(2\lambda)+(\lambda+1.5\lambda+0.5^3)/(3\lambda^2). \label{lem1.5}
\end{align}
When $\lambda\to\infty$, $\eqref{lem1.5}\to0$. Thus, $\lim\sup_{\lambda\to\infty}f(\lambda)\le0.$
\qed

\section{Proof of Lemma~\ref{Lemma 2}.}

We use the following lemmas to prove the positivity of
\begin{align}
    -(x+t+1)\log\Big(1+\frac{s}{1+s}\frac{-2t}{x+2t+1}\Big)-sx\log\Big(1+\frac{1}{1+s}\frac{2t}{x}\Big). \label{lem2}
\end{align}

\begin{lemma}\label{Lemma 2.1}
For any $\alpha>0$, $y\log(1+\frac{\alpha}{y})+\frac{\alpha^2}{2(y+\alpha)}$ is an increasing function of $y>0$.
\end{lemma}

\textbf{Proof of Lemma \ref{Lemma 2.1}.} The differential function is 
\begin{align*}-\log\Big(1-\frac{\alpha}{y+\alpha}\Big)-\frac{\alpha}{y+\alpha}-\frac{1}{2}\Big(\frac{\alpha}{y+\alpha}\Big)^2.
\end{align*}
Because $-\log(1+z)+z-z^2/2$ is a decreasing function, the differential function is positive.

\begin{lemma}\label{Lemma 2.2}
For any $x>0$, $s\in(0,1]$ and $t>0$, $$-\frac{s(1-s)}{x+2t+1}+\frac{(1-s)x+t+1}{(x+\frac{2t}{1+s})(\frac{x+t+1}{s}+\frac{2t}{1+s})}>0.$$
\end{lemma}

\textbf{Proof of Lemma \ref{Lemma 2.2}.} This is equivalent to proving that the following formula is positive:
\begin{align*}
&-s(1-s)\Big(x+\frac{2t}{1+s}\Big)\Big(\frac{x+t+1}{s}+\frac{2t}{1+s}\Big)+((1-s)x+t+1)(x+2t+1)\\
&=-(1-s)x^2-(1-s)\Big(\frac{2t}{1+s}+t+1+\frac{2ts}{1+s}\Big)x-\frac{2ts(1-s)}{1+s}\Big(\frac{t+1}{s}+\frac{2t}{1+s}\Big)\\
&\quad +(1-s)x^2+(t+1+(1-s)(2t+1))x+(t+1)(2t+1)\\
&>\big(-(1-s)(3t+1)+(t+1+(1-s)(2t+1))\big)x-\frac{2ts(1-s)}{1+s}(t+1)\Big(\frac{1}{s}+2\Big)+(t+1)2t\\
&=\big(t+1-(1-s)t\big)x+(t+1)2t\Big(1-\frac{s(1-s)}{1+s}\Big(\frac{1}{s}+2\Big)\Big)>0.
\end{align*}

\begin{lemma}\label{Lemma 2.3}
For any $x>0$, $s\ge1$ and $t>0$, $$\frac{2}{(1+s)^2}\frac{t}{x+2t+1-\frac{2s}{1+s}t}+\frac{2s}{(1+s)^2}\frac{t}{x+t+1+\frac{2t}{1+s}}-\log\Big(1+\frac{2t}{(1+s)(x+t+1)}\Big)>0.$$
\end{lemma}

\textbf{Proof of Lemma \ref{Lemma 2.3}.} Let $y=\frac{t}{x+t+1}$. Then, the lemma is equivalent to $$g(y):=\frac{2}{(1+s)^2}\frac{1}{y^{-1}+\frac{1-s}{1+s}}+\frac{2s}{(1+s)^2}\frac{1}{y^{-1}+\frac{2}{1+s}}-\log\Big(1+\frac{2}{1+s}y\Big)>0.$$
Note that, because $g(0)=0$ and $y\in(0,1)$, we only need to prove that $g'(y)>0$.
In fact,
\begin{align*}
    g'(y)&=\frac{2}{(1+s)^2}\frac{1}{(1+\frac{1-s}{1+s}y)^2}+\frac{2s}{(1+s)^2}\frac{1}{(1+\frac{2}{1+s}y)^2}-\frac{2}{1+s}\frac{1}{1+\frac{2}{1+s}y}\\
    &\ge\frac{2}{(1+s)^2}\Big(\frac{y+1}{1+\frac{2}{1+s}y}\Big)^2+\frac{2s}{(1+s)^2}\frac{1}{(1+\frac{2}{1+s}y)^2}-\frac{2}{1+s}\frac{1}{1+\frac{2}{1+s}y}\\
    &=\frac{2}{(1+s)^2}\frac{1}{(1+\frac{2}{1+s}y)^2}\Big((y+1)^2+s-(1+s)\Big(1+\frac{2}{1+s}y\Big)\Big)>0.
\end{align*}

\noindent
We return to the proof of Lemma \ref{Lemma 2}. We consider the cases $s\le 1$ and $s>1$. For $s\le 1$, we first use Lemma \ref{Lemma 2.1} to deal with the second term of \eqref{lem2}, and then use Lemma \ref{Lemma 2.2} to prove $\eqref{lem2}>0$. For $s>1$, we first use Lemma \ref{Lemma 2.1} to deal with the second term of \eqref{lem2}, and then use Lemma \ref{Lemma 2.3} to prove $\eqref{lem2}>0$. The details of each case are presented below.

\textbf{Case 1:} $s\le1$.

We set $y_1=x$, $y_2=\frac{x+t+1}{s}$, and $\alpha=\frac{2t}{1+s}$. Using Lemma \ref{Lemma 2.1} and $y_2>y_1$, we obtain
$$y_1\log\Big(1+\frac{\alpha}{y_1}\Big)+\frac{\alpha^2}{2(y_1+\alpha)}<y_2\log\Big(1+\frac{\alpha}{y_2}\Big)+\frac{\alpha^2}{2(y_2+\alpha)}.$$
Thus,
\begin{align*}
    &x\log\Big(1+\frac{2t}{(1+s)x}\Big)<-\frac{4t^2}{2(1+s)^2(x+\frac{2t}{1+s})}\\
    &+\frac{x+t+1}{s}\log\Big(1+\frac{2ts}{(x+t+1)(1+s)}\Big)+\frac{4t^2}{2(1+s)^2(\frac{x+t+1}{s}+\frac{2t}{1+s})}.
\end{align*}
Therefore,
\begin{align*}
    &-(x+t+1)\log\Big(1+\frac{s}{1+s}\frac{-2t}{x+2t+1}\Big)-sx\log\Big(1+\frac{1}{1+s}\frac{2t}{x}\Big)\\
    &>-(x+t+1)\log\Big(1+\frac{s}{1+s}\frac{-2t}{x+2t+1}\Big)+s\frac{4t^2}{2(1+s)^2(x+\frac{2t}{1+s})}\\
    &\quad -s\biggl(\frac{x+t+1}{s}\log\Big(1+\frac{2ts}{(x+t+1)(1+s)}\Big)+\frac{4t^2}{2(1+s)^2(\frac{x+t+1}{s}+\frac{2t}{1+s})}\biggr)\\
    &=-(x+t+1)\log\Big(\frac{(x+2t+1-\frac{2ts}{1+s})(x+t+1+\frac{2ts}{1+s})}{(x+t+1)(x+2t+1)}\Big) +\frac{4t^2s}{2(1+s)^2}\Big(\frac{1}{x+\frac{2t}{1+s}}-\frac{1}{\frac{x+t+1}{s}+\frac{2t}{1+s}}\Big)\\
    &=-(x+t+1)\log\Big(1+\frac{\frac{2t^2s}{1+s}-\frac{4t^2s^2}{(1+s)^2}}{(x+t+1)(x+2t+1)}\Big)+\frac{4t^2}{2(1+s)^2}\frac{(1-s)x+t+1}{(x+\frac{2t}{1+s})(\frac{x+t+1}{s}+\frac{2t}{1+s})}\\
    &\ge -(x+t+1)\frac{\frac{2t^2s(1-s)}{(1+s)^2}}{(x+t+1)(x+2t+1)}+\frac{2t^2}{(1+s)^2}\frac{(1-s)x+t+1}{(x+\frac{2t}{1+s})(\frac{x+t+1}{s}+\frac{2t}{1+s})}.
\end{align*}
Based on Lemma \ref{Lemma 2.2}, we know that this value is nonnegative, which completes the proof.

\textbf{Case 2:} $s>1$.

We set $y_1=x$, $y_2=x+t+1$, and $\alpha=\frac{2t}{1+s}$. Using Lemma \ref{Lemma 2.1} and $y_2>y_1$, we obtain
$$y_1\log\Big(1+\frac{\alpha}{y_1}\Big)<y_2\log\Big(1+\frac{\alpha}{y_2}\Big).$$
Thus, $$x\log\Big(1+\frac{2t}{(1+s)x}\Big)<(x+t+1)\log\Big(1+\frac{2t}{(1+s)(x+t+1)}\Big).$$
Therefore, 
\begin{align*}
    &-(x+t+1)\log\Big(1+\frac{s}{1+s}\frac{-2t}{x+2t+1}\Big)-sx\log\Big(1+\frac{1}{1+s}\frac{2t}{x}\Big)\\
    &>-(x+t+1)\log\Big(1+\frac{s}{1+s}\frac{-2t}{x+2t+1}\Big)-s(x+t+1)\log\Big(1+\frac{2t}{(1+s)(x+t+1)}\Big)\\
    &=(x+t+1)\Big(-\log\Big(\frac{x+2t+1-\frac{2s}{1+s}t}{x+2t+1}\Big)-s\log\Big(\frac{x+t+1+\frac{2t}{1+s}}{x+t+1}\Big)\Big) =:(x+t+1)f(s).
\end{align*}
Furthermore, $$f'(s)=\frac{2}{(1+s)^2}\frac{t}{x+2t+1-\frac{2s}{1+s}t}+\frac{2s}{(1+s)^2}\frac{t}{x+t+1+\frac{2t}{1+s}}-\log\Big(1+\frac{2t}{(1+s)(x+t+1)}\Big).$$
From Lemma \ref{Lemma 2.3}, $f'(s)>0.$ We note that $f(1)=0$. Therefore, $f(s)>0$, which completes the proof.

\qed
\end{appendices}

\bibliographystyle{acmtrans.bst}
\bibliography{scholar.bib}
\end{document}